\documentclass[12pt]{article}%
\usepackage{graphicx}
\usepackage[intlimits]{amsmath}
\usepackage{latexsym}
\usepackage{amsfonts}
\usepackage{amssymb}%
\newtheorem{lemma}{Lemma}

\newtheorem{theorem}{Theorem}

\newenvironment{proof}[1][Proof]{\noindent{\textbf {#1}  }}  {\hfill$\Box$\bigskip}

\begin{document}

\title{The minimum number of 4-cliques in a graph with triangle-free complement}
\author{Vladimir Nikiforov\\{\small Department of Mathematical Sciences}\\{\small University of Memphis, Memphis, TN 38152}}
\maketitle

\begin{abstract}
Let $f\left(  n,4,3\right)  $ be minimum number of $4-$cliques in a graph of
order $n$ with independence number $2.$ We show that%
\[
f\left(  n,4,3\right)  \leq\frac{1}{200}n^{4}+O\left(  n^{3}\right)  .
\]
We also show if a graph of order $n$ has independence number $2$ and is close
to regular then it has at least%
\[
\frac{1}{200}n^{4}+o\left(  n^{4}\right)
\]
$4$-cliques.

\end{abstract}

\section{Notation and conventions}

Our notation and terminology are standard (see, e.g. \cite{Bo}). In
particular, all graphs are assumed to be defined on the vertex set $\left[
n\right]  =\left\{  1,2,...n\right\}  .$ For any two adjacent vertices $i$ and
$j$ we write $i\sim j.$ We write $i\nsim j$ if $i$ and $j$ are distinct
nonadjacent vertices. Given a vertex $i,$ $N_{i}$ denotes the set of its
neighbors, $d_{i}$ denotes its degree and $t_{i}$ denotes the number of
triangles containing it.

Given a graph $G$, $t_{3}\left(  G\right)  $ is the number of its triangles;
$t_{4}\left(  G\right)  $ is the number of its $4$-cliques; $t_{3}^{\prime
}\left(  G\right)  $ is the number of all induced subgraphs of order $3$ and
size $2$; $t_{3}^{\prime\prime}\left(  G\right)  $ is the number of all
induced subgraphs of order $3$ and size $1$, and $t_{4}^{\prime}\left(
G\right)  $ is the number of all induced subgraphs of order $4$ and size $5$
(i.e. isomorphic to $K_{4}$ with one edge removed).

\section{Introduction}

Write $f(n,4,3)$ for the function%
\[
f(n,4,3)=\min\{t_{4}\left(  G\right)  :v(G)=n,\text{ }\overline{G}\text{ is
triangle-free}\}.
\]

This function is a particular case of a more general function introduced by
Erd\H{o}s in \cite{Erd1}. In \cite{Ni} we proved that
\[
f(n,4,3)\leq\frac{1}{200}n^{4}+o\left(  n^{4}\right)  .
\]

It happens that this estimate is tight under additional assumptions. We show
that, in fact,%
\[
f(n,4,3)\leq\frac{1}{200}n^{4}+O\left(  n^{3}\right)  .
\]
Moreover, let $G$ be a graph of order $n$ such that $\overline{G}$ is
triangle-free. If $G$ is almost regular, i.e., if
\begin{equation}
\sum_{i\in\left[  n\right]  }\left\vert d_{i}-\frac{2m}{n}\right\vert
=o\left(  n^{2}\right)  ,\label{alreg}%
\end{equation}
then we show that
\[
t_{4}(G)\geq\frac{1}{200}n^{4}+o\left(  n^{4}\right)  .
\]

\section{Main results}

Write $C_{5}\left[  K_{p}\right]  $ for the lexicographic product of the
$5$-cycle with the complete graph on $p$ vertices. Recall that the vertex set
of $G\left[  K_{p}\right]  $ is $v\left(  G\right)  \times v\left(
K_{p}\right)  $ and $(i,x)\sim(j,y)$ iff $i\sim j$ or $x=y$. Observe that the
complement of $C_{5}\left[  K_{p}\right]  $ is triangle-free and
\[
t_{4}\left(  C_{5}\left[  K_{p}\right]  \right)  \geq5\left(  \binom{2p}%
{4}-\binom{p}{4}\right)  =\frac{25}{8}p^{4}-\frac{35}{4}p^{3}+\frac{55}%
{8}p^{2}-\frac{5}{4}p.
\]

Obviously, $f(n,4,3)$ is nondecreasing with respect to $n.$ Let $5p$ be the
smallest multiple of $5$ which is not smaller then $n.$ We have
\begin{align*}
f(n,4,3)  &  \leq5\left(  \binom{2p}{4}-\binom{p}{4}\right)  =\frac{25}%
{8}p^{4}-\frac{35}{4}p^{3}+\frac{55}{8}p^{2}-\frac{5}{4}p\\
&  \leq\frac{25}{8}\left(  \frac{n+4}{5}\right)  ^{4}-\frac{35}{4}\left(
\frac{n+4}{5}\right)  ^{3}+\frac{55}{8}\left(  \frac{n+4}{5}\right)
^{2}-\frac{5}{4}\left(  \frac{n+4}{5}\right) \\
&  =\frac{1}{200}n^{4}+\frac{1}{100}n^{3}-\frac{17}{200}n^{2}-\frac{13}%
{100}n+\frac{1}{5}.
\end{align*}
Our goal to the end of the paper is to prove the following assertion.

\begin{theorem}
\label{th1}If a graph $G$ of order $n$ with independence number $2$ satisfies
(\ref{alreg}) then%
\[
t_{4}\left(  G\right)  \geq\frac{1}{200}n^{4}+o\left(  n^{4}\right)  .
\]

\end{theorem}

The proof consists of the following steps:

\emph{a) }some lemmas establishing properties of graphs with triangle-free complement;

\emph{b) }deduction of a lower bound on $t_{4}\left(  G\right)  $ as a
function of various graph parameters;

\emph{c)} reduction of the bound to a function of vertex degrees, $v\left(
G\right)  $ and $e\left(  G\right)  $;

\emph{d) }replacing the vertex degrees by their mean;

\emph{e) }minimizing the bound with respect to $e\left(  G\right)  $.

\begin{lemma}
\label{le00}For any graph $G=G\left(  n,m\right)  $%
\[
\sum_{i\sim j}d_{i}d_{j}\geq\frac{4m^{3}}{n^{2}}.
\]

\end{lemma}

\begin{proof}
Obviously, we can assume that $G$ contains no isolated vertices. On the one
hand, we have
\[
2\sum_{i\sim j}\frac{1}{\sqrt{d_{i}d_{j}}}\leq\sum_{i\sim j}\frac{1}{d_{i}%
}+\frac{1}{d_{j}}=\sum_{i\in\left[  n\right]  }\frac{d_{i}}{d_{i}}=n.
\]
On the other hand, by the Cauchy-Schwarz inequality,
\[
\left(  \sum_{i\sim j}\sqrt{d_{i}d_{j}}\right)  \left(  \sum_{i\sim j}\frac
{1}{\sqrt{d_{i}d_{j}}}\right)  \geq m^{2}.
\]
Therefore,
\[
\sum_{i\sim j}\sqrt{d_{i}d_{j}}\geq\frac{2m^{2}}{n}.
\]
Hence, by the Cauchy inequality,
\[
\sum_{i\sim j}d_{i}d_{j}\geq\frac{1}{m}\left(  \sum_{i\sim j}\sqrt{d_{i}d_{j}%
}\right)  ^{2}\geq\frac{4m^{3}}{n^{2}}%
\]
completing the proof.
\end{proof}

Let $G=G\left(  n,m\right)  $ be a graph with no independent set on $3$
vertices (i.e., $\overline{G}$ is triangle-free). We shall prove a series of
short lemmas which follow from this assumption.

\begin{lemma}
\label{le0}
\[
6t_{3}\left(  G\right)  =n^{3}-3n^{2}+2n+\sum_{i\in\left[  n\right]  }%
3d_{i}^{2}-3d_{i}n+3d_{i}.
\]

\end{lemma}

\begin{proof}
\ Indeed, this is an expanded form of a well known identity. Observe that
\[
t_{3}^{\prime}\left(  G\right)  +t_{3}^{\prime\prime}\left(  G\right)
=\frac{1}{2}\sum_{i\in\left[  n\right]  }d_{i}\left(  n-1-d_{i}\right)  .
\]
Hence,
\begin{align}
6t_{3}\left(  G\right)   &  =6\binom{n}{3}-6\left(  t_{3}^{\prime}\left(
G\right)  +t_{3}^{\prime\prime}\left(  G\right)  \right) \nonumber\\
&  =n^{3}-3n^{2}+2n-3\sum_{i\in\left[  n\right]  }d_{i}\left(  n-1-d_{i}%
\right) \nonumber\\
&  =n^{3}-3n^{2}+2n+\sum_{i\in\left[  n\right]  }3d_{i}^{2}-3nd_{i}+3d_{i}.
\label{eqt}%
\end{align}

\end{proof}

\begin{lemma}
\label{le1}
\begin{equation}
2t_{3}\left(  G\right)  +t_{3}^{\prime}\left(  G\right)  =\left(  n-2\right)
m-\binom{n}{3}. \label{eq1}%
\end{equation}

\end{lemma}

\begin{proof}
Indeed, trivially,
\[
\binom{n}{3}=t_{3}\left(  G\right)  +t_{3}^{\prime}\left(  G\right)
+t_{3}^{^{\prime\prime}}\left(  G\right)  .
\]
On the other hand,
\[
\left(  n-2\right)  m=3t_{3}\left(  G\right)  +2t_{3}^{\prime}\left(
G\right)  +t_{3}^{^{\prime\prime}}\left(  G\right)  .
\]
Subtracting the last two identities, we obtain (\ref{eq1}).
\end{proof}

\begin{lemma}
\label{le2}
\begin{equation}
8t_{4}\left(  G\right)  +2t_{4}^{\prime}\left(  G\right)  =\sum_{i\in\left[
n\right]  }\left(  d_{i}-2\right)  t_{i}-\sum_{i\in\left[  n\right]  }%
\binom{d_{i}}{3}. \label{eq2}%
\end{equation}

\end{lemma}

\begin{proof}
Applying lemma \ref{le1} to $G\left[  N_{i}\right]  $ for any $i\in\left[
n\right]  ,$ and summing over all vertices, we obtain (\ref{eq2}).
\end{proof}

\begin{lemma}
\label{le3}
\[
t_{4}^{\prime}\left(  G\right)  \leq\sum_{i\nsim j}\binom{\left\vert N_{i}\cap
N_{j}\right\vert }{2}.
\]

\end{lemma}

\begin{proof}
Obviously,
\[
t_{4}^{\prime}\left(  G\right)  =\sum_{i\nsim j}e\left(  \left[  N_{i}\cap
N_{j}\right]  \right)  \leq\sum_{i\nsim j}\binom{\left\vert N_{i}\cap
N_{j}\right\vert }{2}.
\]

\end{proof}

\begin{lemma}
\label{le4}For any two nonadjacent vertices $i$ and $j$
\[
\left|  N_{i}\cap N_{j}\right|  =d_{i}+d_{j}-n+2.
\]

\end{lemma}

\begin{proof}
For every $k,$ if $k\nsim j$ and $k\nsim i$ then $\{i,j,k\}$ is a triangle in
$\overline{G}.$ Therefore,
\[
n-2=\left\vert N_{i}\cup N_{j}\right\vert =d_{i}+d_{j}-\left\vert N_{i}\cap
N_{j}\right\vert .
\]

\end{proof}

\begin{lemma}
\label{le5}
\[
\sum_{i\in\left[  n\right]  }d_{i}t_{i}=\sum_{i\sim j}d_{i}d_{j}-\frac{1}%
{2}\sum_{i\in\left[  n\right]  }d_{i}^{2}-\frac{1}{2}\sum_{i\nsim j}\left(
d_{i}+d_{j}\right)  \left\vert N_{i}\cap N_{j}\right\vert .
\]

\end{lemma}

\begin{proof}
Let $i$ be any vertex. Obviously,
\[
\sum_{j\sim i}d_{j}=d_{i}+2t_{i}+\sum_{j\nsim i}\left\vert N_{i}\cap
N_{j}\right\vert
\]
hence, multiplying both sides by $d_{i},$ we get
\[
d_{i}\sum_{i\sim j}d_{j}=d_{i}^{2}+2d_{i}t_{i}+d_{i}\sum_{j\nsim i}\left\vert
N_{i}\cap N_{j}\right\vert
\]
and summing over all vertices, we obtain
\[
2\sum_{i\sim j}d_{i}d_{j}=\sum_{i\in\left[  n\right]  }d_{i}^{2}+2\sum
_{i\in\left[  n\right]  }d_{i}t_{i}+\sum_{j\nsim i}\left(  d_{i}+d_{j}\right)
\left\vert N_{i}\cap N_{j}\right\vert .
\]

\end{proof}

\begin{proof}
[\textbf{Proof of Theorem \ref{th1}}]By lemma \ref{le2} and lemma \ref{le3},
\begin{align*}
8t_{4}\left(  G\right)   &  =-2t_{4}^{\prime}\left(  G\right)  +\sum
_{i\in\left[  n\right]  }\left(  d_{i}-2\right)  t_{i}-\sum_{i\in\left[
n\right]  }\binom{d_{i}}{3}\\
&  \geq-2\sum_{i\nsim j}\binom{\left\vert N_{i}\cap N_{j}\right\vert }%
{2}-6t_{3}\left(  G\right)  +\sum_{i\in\left[  n\right]  }d_{i}t_{i}%
-\sum_{i\in\left[  n\right]  }\binom{d_{i}}{3}.
\end{align*}
Hence, by lemma \ref{le4} and lemma \ref{le5},
\begin{align*}
8t_{4}\left(  G\right)   &  \geq-\sum_{i\nsim j}\left(  d_{i}+d_{j}%
-n+2\right)  \left(  d_{i}+d_{j}-n+1\right)  -6t_{3}\left(  G\right)
-\sum_{i\in\left[  n\right]  }\binom{d_{i}}{3}\\
&  +\sum_{i\sim j}d_{i}d_{j}-\frac{1}{2}\sum_{i\in\left[  n\right]  }d_{i}%
^{2}-\frac{1}{2}\sum_{i\nsim j}\left(  d_{i}+d_{j}\right)  \left(  d_{i}%
+d_{j}-n+2\right) \\
&  =\sum_{i\sim j}d_{i}d_{j}-\frac{1}{2}\sum_{i\nsim j}\left(  d_{i}%
+d_{j}-n+2\right)  \left(  3d_{i}+3d_{j}-2n+2\right) \\
&  -\frac{1}{6}\sum_{i\in\left[  n\right]  }\left(  d_{i}^{3}-3d_{i}%
^{2}+2d_{i}\right)  -\frac{1}{2}\sum_{i\in\left[  n\right]  }d_{i}^{2}%
-6t_{3}\left(  G\right)  .
\end{align*}
Therefore, we have
\begin{align}
8t_{4}\left(  G\right)   &  \geq\sum_{i\sim j}d_{i}d_{j}-\frac{1}{2}%
\sum_{i\nsim j}\left(  d_{i}+d_{j}-n+2\right)  \left(  3d_{i}+3d_{j}%
-2n+2\right) \label{eq3}\\
&  -\frac{1}{6}\sum_{i\in\left[  n\right]  }\left(  d_{i}^{3}+2d_{i}\right)
-6t_{3}\left(  G\right)  .\nonumber
\end{align}
We also find that
\begin{align*}
&  \sum_{i\nsim j}\left(  d_{i}+d_{j}-n+2\right)  \left(  3d_{i}%
+3d_{j}-2n+2\right) \\
&  =\sum_{i\nsim j}\left(  3d_{i}^{2}+6d_{i}d_{j}-5nd_{i}+8d_{i}+3d_{j}%
^{2}-5nd_{j}+8d_{j}+2n^{2}-6n+4\right) \\
&  =\sum_{i\nsim j}6d_{i}d_{j}+\sum_{i}\left(  3d_{i}^{2}\left(
n-1-d_{i}\right)  -5nd_{i}\left(  n-1-d_{i}\right)  +8d_{i}\left(
n-1-d_{i}\right)  \right)  +\\
&  +\sum_{i\in\left[  n\right]  }\left(  n^{2}-3n+2\right)  \left(
n-1-d_{i}\right) \\
&  =\sum_{i\nsim j}6d_{i}d_{j}+\sum_{i\in\left[  n\right]  }\left(
8nd_{i}^{2}-3d_{i}^{3}-11d_{i}^{2}-6n^{2}d_{i}+16nd_{i}-10d_{i}\right) \\
&  +n^{4}-4n^{3}+5n^{2}-2n.
\end{align*}
Applying this equality to (\ref{eq3}) and afterwards bounding $6t_{3}$ by
lemma \ref{le0}, we get
\begin{align*}
8t_{4}\left(  G\right)   &  \geq\sum_{i\sim j}d_{i}d_{j}-3\sum_{i\nsim j}%
d_{i}d_{j}-\frac{1}{2}\left(  n^{4}-4n^{3}+5n^{2}-2n\right)  -\frac{1}{6}%
\sum_{i}\left(  d_{i}^{3}+2d_{i}\right) \\
&  -6t_{3}+\sum_{i}\left(  \frac{3}{2}d_{i}^{3}-4nd_{i}^{2}+3n^{2}d_{i}%
+\frac{11}{2}d_{i}^{2}-8nd_{i}+5d_{i}\right) \\
&  =\sum_{i\sim j}d_{i}d_{j}-3\sum_{i\nsim j}d_{i}d_{j}-\frac{1}{2}n^{4}%
+n^{3}+\frac{1}{2}n^{2}-n+\\
&  +\sum_{i}\left(  \frac{4}{3}d_{i}^{3}-4nd_{i}^{2}+3n^{2}d_{i}+\frac{5}%
{2}d_{i}^{2}-5nd_{i}+\frac{5}{3}d_{i}\right)
\end{align*}
On the other hand, note that
\begin{align*}
\sum_{i\sim j}d_{i}d_{j}-3\sum_{i\nsim j}d_{i}d_{j}  &  =4\sum_{i\sim j}%
d_{i}d_{j}-\frac{3}{2}\left(  \sum_{i\in\left[  n\right]  }d_{i}\right)
^{2}+\frac{3}{2}\sum_{i\in\left[  n\right]  }d_{i}^{2}\\
&  =4\sum_{i\sim j}d_{i}d_{j}-6m^{2}+\frac{3}{2}\sum_{i\in\left[  n\right]
}d_{i}^{2}%
\end{align*}
Thus,
\begin{align*}
8t_{4}\left(  G\right)   &  \geq4\sum_{i\sim j}d_{i}d_{j}-6m^{2}-\frac{1}%
{2}n^{4}+n^{3}+\frac{1}{2}n^{2}-n+\\
&  +\sum_{i\in\left[  n\right]  }\left(  \frac{4}{3}d_{i}^{3}-4nd_{i}%
^{2}+3n^{2}d_{i}+4d_{i}^{2}-5nd_{i}+\frac{5}{3}d_{i}\right)
\end{align*}
Dropping the low order terms, we see that
\[
8t_{4}\geq4\sum_{i\sim j}d_{i}d_{j}-6m^{2}-\frac{1}{2}n^{4}+\sum_{i\in\left[
n\right]  }\left(  \frac{4}{3}d_{i}^{3}-4nd_{i}^{2}+3n^{2}d_{i}\right)
+O\left(  n^{3}\right)  .
\]
Due to (\ref{alreg}), we find that
\begin{align*}
8t_{4}  &  \geq4m\frac{4m^{2}}{n^{2}}-6m^{2}-\frac{1}{2}n^{4}+\frac{4}{3}%
\frac{8m^{3}}{n^{2}}-4n\frac{4m^{2}}{n}+6n^{2}m+O\left(  n^{3}\right) \\
&  =\frac{80}{3}\frac{m^{3}}{n^{2}}-22m^{2}+6mn^{2}-\frac{1}{2}n^{4}+O\left(
n^{3}\right)  .
\end{align*}
Since the expression
\[
\frac{80}{3}\frac{m^{3}}{n^{2}}-22m^{2}+6mn^{2}%
\]
attains its minimum at
\[
m=\frac{3}{10}n^{2}+o\left(  n^{2}\right)  ,
\]
the desired inequality follows.
\end{proof}

\textbf{Acknowledgment.} The author is grateful to Andrew Thomason for his
encouraging attention.

\end{document}